\documentclass{article}
\usepackage[utf8]{inputenc}
\usepackage{amsmath}
\usepackage{tikz}
\usepackage{graphicx}

\title{Optimum Strategies for the Game Super Six}
\author{R\"udiger Jehn }
\date{September 2021}

\usepackage[colorlinks,citecolor=blue,urlcolor=blue,bookmarks=false,hypertexnames=true]{hyperref} 
\begin{document}

\maketitle

{\bf Abstract}\\

In the game "Super Six", after successfully getting rid of a stick by rolling with the die a number that is not occupied, the player has the choice to continue to roll the die or to stop and to hand over the die to their opponent.
The probability to win depends on the situation and on the chosen strategy. In this paper, the optimum strategy to maximise the winning probability is derived. If there are 1 or 2 sticks on the lid, it is always best to continue to roll the die. If there are 5 sticks on the lid it is best to stop. If there are 4 sticks on the lid it is best to stop unless both players have one stick in their hands. If there are 3~sticks on the lid, the optimum strategy depends on the number of sticks in the hands of the two players.

\section{Introduction}
The rules of Super Six \cite{super6} for two players are simple. Game equipment consists of a six-sided die, a number of sticks, and a box whose lid has six holes. The holes numbered 1 through 5 are shallow, and a stick placed in any one of them will stand up in it; hole $\#$6 goes all the way through the lid so that any stick placed in it falls into the box and is out of play. Initially, an even number of sticks are divided evenly between the two players. The goal is to get rid of all one's sticks before the other player does.\\

The players take turns. On each turn, the player whose turn it is rolls the die and places a stick in the numbered hole that matches the number on the die (e.g., a player who rolls a 4 then places a stick in hole $\#$4). The player may roll and place a stick for each roll as many times as desired until rolling a number that is already filled by a stick. When this occurs, the player must take that stick in hand, and play passes to the opposing player.\\

The game proceeds with players taking turns and ends when one player has run out of sticks. The only freedom that the players have is the decision whether to continue rolling the die or not after successfully placing a stick.\\

In this paper we will calculate the winning probabilities for all possible strategies. Based on these probabilities we will search for the best strategy in all situations where the player has the choice to keep rolling the die or to stop.

\section{Strategies}

A game situation can be described by 3 numbers: $j$, $k$ and $l$ where $j$ is the number of sticks on the lid, $k$ is the number of sticks in the hand of player 1 and $l$ is the number of sticks in the hand of the opponent, player 2. Let $i=j+k+l$ be the total number of sticks in the current game. We will use the notation ($j/k/l$) to describe a situation.\\

A strategy $S_i$ consists of $n$ sub-strategies: ($s_n s_{n-1} ... s_2 s_1$) where each sub-strategy is either 0 or 1, 0 meaning that the player will stop rolling the die at a given situation and 1 meaning that the player will continue.
$P_A^{S_i}$ is the probability to win at situation $A$ when playing the strategy $S_i$.\\

We assume that a player will always continue to roll the die if the lid is empty ($j=0$).\\

The situations (1/k/1) do not require a strategy, because a player has only a choice when he had successfully placed a stick. If he had rolled a 6, the previous situation was (1/k+1/1) and we are back to the same kind of situations. If he had rolled any other number, the previous situation was (0/k+1/1). But to reach this situation, the only possible previous situation is (0/k+2/1) given the assumption that the opponent does not stop rolling the die if the lid is empty. Again, the chain of possible previous situations does not lead to any starting position (0/k/k) where $k>1$.\\

This means we only have 2 situations where a decision is required when there are 4 sticks in the game and we have four possible strategies $S_4$: (00), (01), (10) and (11). There are 5 sub-strategies when there are 5 sticks, 9 sub-strategies when there are 6 sticks and 5 more sub-strategies with every extra stick in the game. With 7 sticks in the game, the player can choose for instance the strategy $S_7$ = (0101/1010/110/10/1) where the slashes are added to increase the readability. The 1 at the right means continue rolling the die at situation (5/1/1), the 10 means continue at situation (4/1/2) and stop at (4/2/1) and so forth. The slashes separate the sub-strategies according to the number of sticks on the lid.

\section{Three and Four sticks in the game}

If there are 3 sticks left in the game (i.e., held by a player or standing on the lid), the winning probabilities are:\\
$$P_{(0/1/2)} = 1$$ 
$$P_{(0/2/1)} = 31/36$$ 
$$P_{(1/1/1)} = 5/6$$ 

If there are 4 sticks left in the game, there are two situations where a decision may be required. Position A is when the player arrives at (1/1/2) and position B when they arrive at (2/1/1). It is simple to demonstrate that in both situations it is best to continue rolling the die.\\

$P_A$, the probability to win at (1/1/2), is larger than $\frac{5}{6}$, because in 5 out of 6 cases the player wins immediately and in the sixth case they still have a chance to win, because the opponent still has to win their game at a situation of (0/2/2).\\ 

But if the player decides to stop their opponent will be in the situation C=(1/2/1) for which the probability $P_C$ to win can be easily approximated by analyzing the different possibilities: 

\begin{enumerate}
    \item the opponent rolls a 6 and arrives at (1/1/1) with a winning probability of $\frac{5}{6}$
    \item the opponent rolls one of the 4 numbers that are free and arrives at (2/1/1) with a winning probability of $> \frac{4}{6}$
    \item the opponent rolls the only number that is occupied on the lid and loses immediately
\end{enumerate}

Therefore $P_C$ is larger than $\frac{1}{6} \cdot \frac{5}{6} + \frac{4}{6} \cdot \frac{4}{6} = \frac{21}{36}$. If the player stops their chances drop from more than $\frac{5}{6}$ to less than $\frac{15}{36}$.\\

If the player arrives at (2/1/1) and has to take a decision, the situation is also obvious. They have an immediate chance of winning of $\frac{4}{6}$ and if they hand over the die to the opponent, they also hand over the winning probability of more than $\frac{4}{6}$ to the opponent.\\

In both situations it is optimal to continue to roll the die. The winning probability at position (1/1/2) depends on the winning probability at position (0/2/2) and vice versa because you can go from both positions to the other position. Therefore we need to solve a system of equations. 
D is the situation (0/2/2). In this situation two cases can happen:\\

$\frac{1}{6}$: arriving at situation (0/1/2) and winning immediately\\

$\frac{5}{6}$: arriving at situation A = (1/1/2)\\

Also in situation A = (1/1/2) two cases can happen.\\ 

$\frac{5}{6}$: winning immediately\\

$\frac{1}{6}$: arriving at situation D = (0/2/2) but swapping the die \\

We can see that in one out of six cases the player in situation A will have to pick up the single stick on the board and we return to the situation D but with the difference that now it is the turn of the other player. We get the following equations:

$$P_D = \frac{1}{6} + \frac{5}{6} P_A$$

$$P_A = \frac{5}{6} + \frac{1}{6} \cdot (1 – P_D)$$ 

$1 – P_D$ is the probability for the player to win when the opponent is given the die at the situation D = (0/2/2).\\ 

Hence we have to solve this system of equations:
$$
\begin{pmatrix}
6 & -5 \\
1 & 6  
\end{pmatrix}
\begin{pmatrix}
P_D\\
P_A  
\end{pmatrix}
=
\begin{pmatrix}
1\\
6  
\end{pmatrix}
$$
which yields $P_D$ = $\frac{36}{41}$ and $P_A$ = $\frac{35}{41}$. The winning probability at situation (0/3/1) can be easily calculated, but this position will never be reached unless a player stops rolling the die when there are zero sticks on the lid. The other two probabilities are:

$$P_B = \frac{4}{6} + \frac{2}{6} \cdot (1 – P_A) = \frac{88}{123}$$ 

$$P_C = \frac{1}{6} P_{(1/1/1)} + \frac{4}{6} P_B = 0.616$$ 

Figure \ref{fig:4sticks} summarizes the winning probabilities. 

\begin{figure}[h]
\begin{center}
\begin{tikzpicture}

\draw[blue, very thick] (0,0) rectangle (2,1);
\draw[blue, very thick] (2,0) rectangle (4,1);
\draw[blue, very thick] (4,0) rectangle (6,1);
\draw[blue, very thick] (1,1) rectangle (3,2);
\draw[blue, very thick] (3,1) rectangle (5,2);
\draw[blue, very thick] (2,2) rectangle (4,3);

\node (c) at (1,0.7) {(0.657)};
\node (c) at (0.5,0.3) {\small 3};
\node (c) at (1.5,0.3) {\small 1};

\node (c) at (3,0.7) {0.878};
\node (c) at (3.8,0.8) {D};
\node (c) at (2.5,0.3) {\small 2};
\node (c) at (3.5,0.3) {\small 2};

\node (c) at (5,0.7) {1.000};
\node (c) at (4.5,0.3) {\small 1};
\node (c) at (5.5,0.3) {\small 3};
\node (c) at (-1,0.3) {\small 0}; 

\node (c) at (2,1.7) {0.616};
\node (c) at (2.8,1.8) {C};
\node (c) at (1.5,1.3) {\small 2};
\node (c) at (2.5,1.3) {\small 1};

\node (c) at (4,1.7) {0.853};
\node (c) at (4.8,1.8) {A};
\node (c) at (3.5,1.3) {\small 1};
\node (c) at (4.5,1.3) {\small 2};
\node (c) at (-1,1.3) {\small 1};  

\node (c) at (3,2.7) {0.715};
\node (c) at (3.8,2.8) {B};
\node (c) at (2.5,2.3) {\small 1};
\node (c) at (3.5,2.3) {\small 1};
\node (c) at (-1,2.3) {\small 2};  

\end{tikzpicture}
\end{center}
\vspace{-5mm}
\caption{The winning probabilities in the game "Super Six" when there are 4 sticks in the game. The numbers on the left indicate the number of sticks on the lid, the two numbers below the probabilities indicate how many sticks the two players have in hand.}
\label{fig:4sticks}
\end{figure}

\section{Five sticks in the game}

With 5 sticks in the game, there are 5 situations where a player may be required to take a decision. When arriving at the situation (3/1/1) the strategy is obvious:\\

In 3 out of 6 cases the player who rolls the die wins immediately and in the other 3 cases they still have a chance to win, because the opponent still has to win his game at a situation of (2/1/2). So the chances to win are larger than 50 $\%$ and that advantage should not be given to the opponent.\\ 

When a player arrives at (2/1/2), the probability to win is larger than $\frac{4}{6}$, because in 4 out of 6 cases they win immediately and in the other case they still have a chance to win.\\ 

Handing over the die to the opponent would put them in situation (2/2/1) with 3 possible outcomes.\\ 

$\frac{1}{6}$: arriving at (2/1/1) with a winning probability of 88/123  (see 4 sticks) \\

$\frac{3}{6}$: arriving at (3/1/1) with a winning probability $> 0.5$ (see first situation) \\

$\frac{2}{6}$: arriving at (1/3/1) with a winning probability $> 0$ \\

We can approximate: $P_{(2/2/1)} > \frac{1}{6} \cdot \frac{88}{123} + 0.5 \cdot 0.5 = 0.369 > \frac{2}{6}$. Handing over the die to the opponent will raise their chances from below $\frac{2}{6}$ to more than $\frac{2}{6}$ and therefore the right strategy at (2/1/2) is to continue the game.\\

Also for the other 3 situations it is pretty simple to demonstrate that to continue rolling the die is always the best strategy. The proofs are left as a small exercise.\\

Having sorted out the right strategies, the winning probabilities for each situation can now be calculated. There are 5~situations through which the game can pass repeatedly. Therefore we get a system of 5~equations: 

$$ P_{(0/3/2)} = \frac{1}{6} \cdot P_{(0/2/2)} + \frac{5}{6} \cdot P_{(1/2/2)} $$
$$ P_{(1/2/2)} = \frac{1}{6} \cdot P_{(1/1/2)} + \frac{4}{6} \cdot P_{(2/1/2)} + \frac{1}{6} \cdot (1 - P_{(0/2/3)}) $$
$$ P_{(2/1/2)} = \frac{4}{6} + \frac{2}{6} \cdot (1 - P_{(1/2/2)}) $$
$$ P_{(0/2/3)} = \frac{1}{6} + \frac{5}{6} \cdot P_{(1/1/3)} $$
$$ P_{(1/1/3)} = \frac{5}{6} + \frac{1}{6} \cdot (1 - P_{(0/3/2)}) $$

If we write the relations of these 5 situations in matrix form, we have to solve this system of equations: 
$$
\begin{pmatrix}
6 & -5 & 0 & 0 & 0\\
0 & 6 & -4 & 1 & 0\\
0 & 2 & 6 & 0 & 0\\
0 & 0 & 0 & 6 & -5\\
1 & 0 & 0 & 0 & 6
\end{pmatrix}
\begin{pmatrix}
P_{(0/3/2)}\\
P_{(1/2/2)}\\
P_{(2/1/2)}\\
P_{(0/2/3)}\\
P_{(1/1/3)}  
\end{pmatrix}
=
\begin{pmatrix}
\frac{36}{41}\\
\frac{76}{41}\\
6\\
1\\
6  
\end{pmatrix}
$$
The solution vector is
$$(45324,   43164,   49531,   57624,   56365) / 63919$$

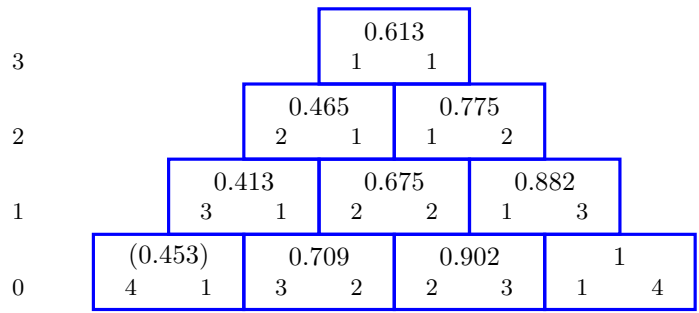
\begin{figure}[h]
\begin{center}
\begin{tikzpicture}

\draw[blue, very thick] (0,0) rectangle (2,1);
\draw[blue, very thick] (2,0) rectangle (4,1);
\draw[blue, very thick] (4,0) rectangle (6,1);
\draw[blue, very thick] (6,0) rectangle (8,1);
\draw[blue, very thick] (1,1) rectangle (3,2);
\draw[blue, very thick] (3,1) rectangle (5,2);
\draw[blue, very thick] (5,1) rectangle (7,2);
\draw[blue, very thick] (2,2) rectangle (4,3);
\draw[blue, very thick] (4,2) rectangle (6,3);
\draw[blue, very thick] (3,3) rectangle (5,4);
\node (c) at (1,0.7) {(0.453)};
\node (c) at (0.5,0.3) {\small 4};
\node (c) at (1.5,0.3) {\small 1};

\node (c) at (3,0.7) {0.709};
\node (c) at (2.5,0.3) {\small 3};
\node (c) at (3.5,0.3) {\small 2};

\node (c) at (5,0.7) {0.902};
\node (c) at (4.5,0.3) {\small 2};
\node (c) at (5.5,0.3) {\small 3};

\node (c) at (7,0.7) {1};
\node (c) at (6.5,0.3) {\small 1};
\node (c) at (7.5,0.3) {\small 4};
\node (c) at (-1,0.3) {\small 0}; 

\node (c) at (2,1.7) {0.413};
\node (c) at (1.5,1.3) {\small 3};
\node (c) at (2.5,1.3) {\small 1};

\node (c) at (4,1.7) {0.675};
\node (c) at (3.5,1.3) {\small 2};
\node (c) at (4.5,1.3) {\small 2};

\node (c) at (6,1.7) {0.882};
\node (c) at (5.5,1.3) {\small 1};
\node (c) at (6.5,1.3) {\small 3};
\node (c) at (-1,1.3) {\small 1}; 

\node (c) at (3,2.7) {0.465};
\node (c) at (2.5,2.3) {\small 2};
\node (c) at (3.5,2.3) {\small 1};

\node (c) at (5,2.7) {0.775};
\node (c) at (4.5,2.3) {\small 1};
\node (c) at (5.5,2.3) {\small 2};
\node (c) at (-1,2.3) {\small 2};  

\node (c) at (4,3.7) {0.613};
\node (c) at (3.5,3.3) {\small 1};
\node (c) at (4.5,3.3) {\small 1};
\node (c) at (-1,3.3) {\small 3};  

\end{tikzpicture}
\end{center}
\vspace{-5mm}
\caption{The winning probabilities when there are 5 sticks in the game. The numbers on the left indicate the number of sticks on the lid, the two numbers below the probabilities indicate how many sticks the two players have in hand.}
\label{fig:5sticks}
\end{figure}

All winning probabilities are listed in Figure \ref{fig:5sticks}.\\

\section{Six sticks in the game}

It is still quite simple to prove that it is best to continue to roll the die if 1 or 2 sticks are on the lid by approximating the winning probabilities, but when there are 3 sticks on the lid it becomes very cumbersome. Instead we calculate the winning probabilities for all strategies and will take our conclusions from there.\\

We must first clarify if it is possible that the optimum strategy of player 1 depends on the strategy of player 2. Let us for a moment assume that the best strategy for player 1 is to continue at a situation A if player 2 has the strategy to stop there and vice versa. But how can player 2 select the strategy to stop at situation A if there was no independent optimum strategy? Player 2 is in the same situation as player 1 and could not choose a best strategy. Therefore, the choice cannot depend on the opponent's strategy and an optimum strategy exists in each situation.\\

With 6 sticks in the game, there are 15 possible situations. (0/1/5) is trivial, because the game ends with the next move. (0/5/1) is a position which cannot be reached unless a player stops rolling the die when the lid is empty. 13 situations remain through which a game can iterate for ever as long as no 6 is rolled. In 9 situations a decision may be required. As an example we start with situation (0/4/2). The winning probability depends if the player continues after rolling any number other than a 6.\\

If they continue, $P_{(0/4/2)} = \frac{1}{6} + \frac{5}{6} \cdot P_{(1/3/2)}$\\

If they stop, $P_{(0/4/2)} = \frac{1}{6} + \frac{5}{6} \cdot (1 - P_{(1/2/3)})$\\

For the situation (1/4/1) we get these two possibilities (note when they hit the only occupied number, they have lost immediately, because the opponent has only 1 stick and an empty board):\\

If they continue, $P_{(1/4/1)} = \frac{1}{6} \cdot P_{(1/3/1)} + \frac{4}{6} \cdot P_{(2/3/1)}$\\

If they stop, $P_{(1/4/1)} = \frac{1}{6} \cdot P_{(1/3/1)} + \frac{4}{6} \cdot (1 - P_{(2/1/3)})$\\

For situation (4/1/1) there is no choice, the winning probability is always
$$P_{(4/1/1)} = \frac{2}{6} + \frac{4}{6} \cdot (1 - P_{(3/1/2)})$$

The corresponding equations need to be established for the other 10 situations as well. For each of the $2^9 = 512$ strategies, we get a different system of 13 equations, which are solved with a computer program.\\

The next step is to compare the winning probabilities $P_{(0/4/2)}^{1\bar{S}}$ if the player continues, with the winning probabilities $P_{(0/4/2)}^{0\bar{S}}$ if the player stops when the player rolls any number other than a 6 at the position (0/4/2). This comparison is made for all 256 strategies $\bar{S}$ at the other 8 decision points. It turns out that in each of the 256 comparisons $P_{(0/4/2)}^{1\bar{S}} > P_{(0/4/2)}^{0\bar{S}}$ which means that continuing to roll the die in this situation is the best strategy (no wonder, there are 5 out of 6 chances to get rid of another stick).\\

These 256 comparisons are made for each of the 9 situations where a strategy is required. It turns out that in all 9 situations, the best strategy is always to continue rolling the die, whenever a decision is to be taken. Surprisingly even at the situation (4/1/1) it is best to continue, even if the immediate chances to get rid of a stick is only 2 out of 6. The winning probability is 0.524 when the player continues and drops to 0.476 if they hand over the die to the opponent.\\

\begin{figure}[h]
\begin{center}
\begin{tikzpicture}

\draw[blue, very thick] (0,0) rectangle (2,1);
\draw[blue, very thick] (2,0) rectangle (4,1);
\draw[blue, very thick] (4,0) rectangle (6,1);
\draw[blue, very thick] (6,0) rectangle (8,1);
\draw[blue, very thick] (8,0) rectangle (10,1);
\filldraw[fill=red!8] (3,1) rectangle (9,2);
\draw[blue, very thick] (1,1) rectangle (3,2);
\draw[blue, very thick] (3,1) rectangle (5,2);
\draw[blue, very thick] (5,1) rectangle (7,2);
\draw[blue, very thick] (7,1) rectangle (9,2);
\filldraw[fill=red!8] (2,2) rectangle (8,3);
\draw[blue, very thick] (2,2) rectangle (4,3);
\draw[blue, very thick] (4,2) rectangle (6,3);
\draw[blue, very thick] (6,2) rectangle (8,3);
\filldraw[fill=red!8] (3,3) rectangle (7,4);
\draw[blue, very thick] (3,3) rectangle (5,4);
\draw[blue, very thick] (5,3) rectangle (7,4);
\filldraw[fill=red!8] (4,4) rectangle (6,5);
\draw[blue, very thick] (4,4) rectangle (6,5);
\node (c) at (1,0.7) {(0.293)};
\node (c) at (0.5,0.3) {\small 5};
\node (c) at (1.5,0.3) {\small 1};

\node (c) at (3,0.7) {0.541};
\node (c) at (2.5,0.3) {\small 4};
\node (c) at (3.5,0.3) {\small 2};

\node (c) at (5,0.7) {0.767};
\node (c) at (4.5,0.3) {\small 3};
\node (c) at (5.5,0.3) {\small 3};

\node (c) at (7,0.7) {0.925};
\node (c) at (6.5,0.3) {\small 2};
\node (c) at (7.5,0.3) {\small 4};

\node (c) at (9,0.7) {1};
\node (c) at (8.5,0.3) {\small 1};
\node (c) at (9.5,0.3) {\small 5};
\node (c) at (-1,0.3) {\small 0}; 

\node (c) at (2,1.7) {0.261};
\node (c) at (1.5,1.3) {\small 4};
\node (c) at (2.5,1.3) {\small 1};

\node (c) at (4,1.7) {0.507};
\node (c) at (3.5,1.3) {\small 3};
\node (c) at (4.5,1.3) {\small 2};

\node (c) at (6,1.7) {0.740};
\node (c) at (5.5,1.3) {\small 2};
\node (c) at (6.5,1.3) {\small 3};

\node (c) at (8,1.7) {0.910};
\node (c) at (7.5,1.3) {\small 1};
\node (c) at (8.5,1.3) {\small 4};
\node (c) at (-1,1.3) {\small 1};

\node (c) at (3,2.7) {0.288};
\node (c) at (2.5,2.3) {\small 3};
\node (c) at (3.5,2.3) {\small 1};

\node (c) at (5,2.7) {0.573};
\node (c) at (4.5,2.3) {\small 2};
\node (c) at (5.5,2.3) {\small 2};

\node (c) at (7,2.7) {0.831};
\node (c) at (6.5,2.3) {\small 1};
\node (c) at (7.5,2.3) {\small 3};
\node (c) at (-1,2.3) {\small 2};  

\node (c) at (4,3.7) {0.361};
\node (c) at (3.5,3.3) {\small 2};
\node (c) at (4.5,3.3) {\small 1};

\node (c) at (6,3.7) {0.714};
\node (c) at (5.5,3.3) {\small 1};
\node (c) at (6.5,3.3) {\small 2};
\node (c) at (-1,3.3) {\small 3};  

\node (c) at (5,4.7) {0.524};
\node (c) at (4.5,4.3) {\small 1};
\node (c) at (5.5,4.3) {\small 1};
\node (c) at (-1,4.3) {\small 4};  

\end{tikzpicture}
\end{center}
\vspace{-5mm}
\caption{The winning probabilities when there are 6 sticks in the game. The numbers on the left indicate the number of sticks on the lid, the two numbers below the probabilities indicate how many sticks the two players have in hand. In the 9 shaded situations a decision whether to continue or to stop may be required.}
\label{fig:6sticks}
\end{figure}
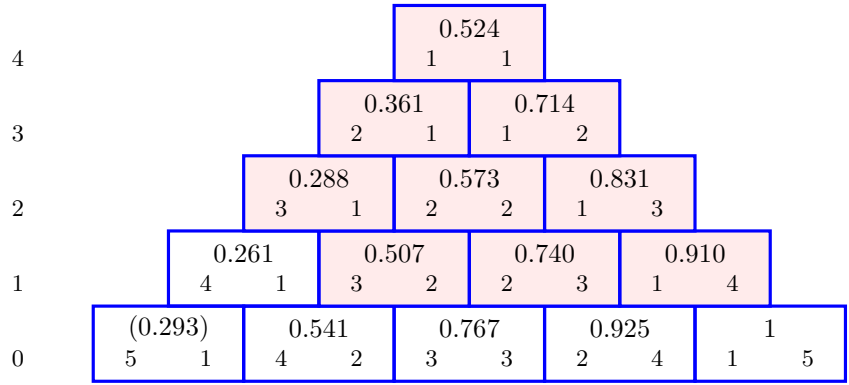

The pyramid with the winning probabilities is illustrated in Figure \ref{fig:6sticks}.

\section{Seven sticks in the game}

The same procedure as with 6 sticks needs to be repeated. However, in contrast to games with 6 or less sticks it can happen that the optimum strategy at a given position is not always superior if at other positions inferior strategies are chosen. For instance:\\

$P_{(4,1,2)}^{(0101/1010/110/10/1)} = 0.1759$\\ 

$P_{(4,1,2)}^{(0101/1010/{\bf 0}10/10/1)} = 0.1765$\\ 

If both players had chosen to play the clearly sub-optimal strategy (0101/ 1010/x10/10/1) at the other decision points, then indeed when you go from (4,1,2) to (3,1,3) it would be slightly better to stop rolling the die to increase your winning chances from 0.1759 to 0.1765.\\

Therefore, the optimum strategies are first determined for situations where 1, 2, 4 or 5 sticks are on the lid. For example, if only 1 stick is on the lid, it turns out by comparing the winning probabilities that, no matter what strategy is chosen at all other situations, the best strategy in these situation with 1 stick on the lid is to continue to roll the die. The same holds true for situations with 2 sticks on the lid: the optimum strategy is to continue to roll the die.\\

In case of 4 or 5 sticks on the lid, the comparison of the winning probabilities, reveals that stop rolling the die is always the best strategy, no matter what strategy is chosen at the other situations.\\

Finally the situations with 3 sticks on the lid need to be analysed. But for this analysis, we can assume that at the other situation with 1, 2, 4 or 5 sticks on the lid the optimum strategies are applied. With 7 sticks in the game and 3 sticks on the lid, there are 3 possible decision points ((3,1,3), (3,2,2) and (3,3,1)) and hence 8 possible strategies. We start comparing the winning probabilities at (2,2,3) for all 4 pairs: Strategy (100) wins over (000), (101) wins over (001), (110) wins over (010) and (111) wins over (011). Hence the best strategy when arriving from (2,2,3) to (3,1,3) is to continue to roll the die.\\

Next we compare the 4 possible strategies, when the player arrives from (2,3,2) to (3,2,2) assuming that we continue rolling the die at situation (3,1,3). Strategy (110) wins over (100) and (111) wins over (101). Hence the best strategy in this situation is to continue to roll the die. Finally we need to compare strategy (111) with (110) when arriving from (2,4,1) to (3,3,1). The first is better than the second and therefore also in this situation it is best to continue.\\

To summarize, with 7 sticks in the game, the best strategy is to continue to roll the die if 1 to 3 sticks are on the lid and to stop if there are 4 or 5 sticks on the lid. This optimum strategy is denoted as
$$\hat{S}_7 = (1111/1111/111/00/0).$$

\begin{figure}[h]
\begin{center}
\begin{tikzpicture}

\draw[blue, very thick] (0,0) rectangle (2,1);
\draw[blue, very thick] (2,0) rectangle (4,1);
\draw[blue, very thick] (4,0) rectangle (6,1);
\draw[blue, very thick] (6,0) rectangle (8,1);
\draw[blue, very thick] (8,0) rectangle (10,1);
\draw[blue, very thick] (10,0) rectangle (12,1);
\filldraw[fill=red!8] (3,1) rectangle (11,2);
\draw[blue, very thick] (1,1) rectangle (3,2);
\draw[blue, very thick] (3,1) rectangle (5,2);
\draw[blue, very thick] (5,1) rectangle (7,2);
\draw[blue, very thick] (7,1) rectangle (9,2);
\draw[blue, very thick] (9,1) rectangle (11,2);
\filldraw[fill=red!8] (2,2) rectangle (10,3);
\draw[blue, very thick] (2,2) rectangle (4,3);
\draw[blue, very thick] (4,2) rectangle (6,3);
\draw[blue, very thick] (6,2) rectangle (8,3);
\draw[blue, very thick] (8,2) rectangle (10,3);
\filldraw[fill=red!8] (3,3) rectangle (9,4);
\draw[blue, very thick] (3,3) rectangle (5,4);
\draw[blue, very thick] (5,3) rectangle (7,4);
\draw[blue, very thick] (7,3) rectangle (9,4);
\filldraw[fill=red!8] (4,4) rectangle (8,5);
\draw[blue, very thick] (4,4) rectangle (6,5);
\draw[blue, very thick] (6,4) rectangle (8,5);
\filldraw[fill=red!8] (5,5) rectangle (7,6);
\draw[blue, very thick] (5,5) rectangle (7,6);
\node (c) at (1,0.7) {(0.190)};
\node (c) at (0.5,0.3) {\small 6};
\node (c) at (1.5,0.3) {\small 1};

\node (c) at (3,0.7) {0.401};
\node (c) at (2.5,0.3) {\small 5};
\node (c) at (3.5,0.3) {\small 2};

\node (c) at (5,0.7) {0.627};
\node (c) at (4.5,0.3) {\small 4};
\node (c) at (5.5,0.3) {\small 3};

\node (c) at (7,0.7) {0.819};
\node (c) at (6.5,0.3) {\small 3};
\node (c) at (7.5,0.3) {\small 4};

\node (c) at (9,0.7) {0.944};
\node (c) at (8.5,0.3) {\small 2};
\node (c) at (9.5,0.3) {\small 5};

\node (c) at (11,0.7) {1.000};
\node (c) at (10.5,0.3) {\small 1};
\node (c) at (11.5,0.3) {\small 6};

\node (c) at (2,1.7) {0.169};
\node (c) at (1.5,1.3) {\small 5};
\node (c) at (2.5,1.3) {\small 1};

\node (c) at (4,1.7) {0.373};
\node (c) at (3.5,1.3) {\small 4};
\node (c) at (4.5,1.3) {\small 2};

\node (c) at (6,1.7) {0.599};
\node (c) at (5.5,1.3) {\small 3};
\node (c) at (6.5,1.3) {\small 3};

\node (c) at (8,1.7) {0.798};
\node (c) at (7.5,1.3) {\small 2};
\node (c) at (8.5,1.3) {\small 4};

\node (c) at (10,1.7) {0.933};
\node (c) at (9.5,1.3) {\small 1};
\node (c) at (10.5,1.3) {\small 5};

\node (c) at (3,2.7) {0.188};
\node (c) at (2.5,2.3) {\small 4};
\node (c) at (3.5,2.3) {\small 1};

\node (c) at (5,2.7) {0.419};
\node (c) at (4.5,2.3) {\small 3};
\node (c) at (5.5,2.3) {\small 2};

\node (c) at (7,2.7) {0.668};
\node (c) at (6.5,2.3) {\small 2};
\node (c) at (7.5,2.3) {\small 3};

\node (c) at (9,2.7) {0.876};
\node (c) at (8.5,2.3) {\small 1};
\node (c) at (9.5,2.3) {\small 4};

\node (c) at (4,3.7) {0.236};
\node (c) at (3.5,3.3) {\small 3};
\node (c) at (4.5,3.3) {\small 1};

\node (c) at (6,3.7) {0.512};
\node (c) at (5.5,3.3) {\small 2};
\node (c) at (6.5,3.3) {\small 2};

\node (c) at (8,3.7) {0.790};
\node (c) at (7.5,3.3) {\small 1};
\node (c) at (8.5,3.3) {\small 3};

\node (c) at (5,4.7) {0.319};
\node (c) at (4.5,4.3) {\small 2};
\node (c) at (5.5,4.3) {\small 1};

\node (c) at (7,4.7) {0.658};
\node (c) at (6.5,4.3) {\small 1};
\node (c) at (7.5,4.3) {\small 2};

\node (c) at (6,5.7) {0.451};
\node (c) at (5.5,5.3) {\small 1};
\node (c) at (6.5,5.3) {\small 1};

\end{tikzpicture}
\end{center}
\vspace{-5mm}
\caption{The winning probabilities when there are 7 sticks in the game. The two numbers below the probabilities indicate how many sticks the two players have in hand. In the 14 shaded situations a decision whether to continue or to stop may be required.}
\label{fig:7sticks}
\end{figure}
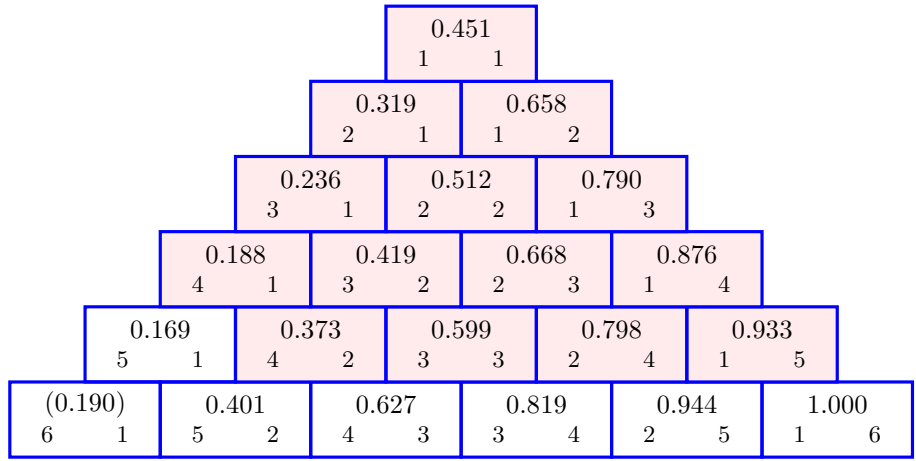

The pyramid with the winning probabilities is illustrated in Figure \ref{fig:7sticks}.

\section{Eight to Twelve sticks in the game}

The same procedure is now repeated for more sticks in the game. First it is shown that with 1 and 2 sticks on the lid, it is best to continue to roll the die, then with 4 and 5 sticks, it is best to stop and the last step is to prove that also for 3 sticks on the lid it is always best to continue to roll the die. The optimum strategies are:\\ 

$\hat{S}_{8} = (11111/11111/1111/000/00)$

$\hat{S}_{9} = (111111/111111/11111/0000/000)$

$\hat{S}_{10} = (1111111/1111111/111111/00000/0000)$

$\hat{S}_{11} = (11111111/11111111/1111111/000000/00000)$

$\hat{S}_{12} = (111111111/111111111/11111111/0000000/000000)$

\section{Thirteen and more sticks in the game}

With 13 sticks in the game, it turns out that it is not always the best strategy to continue to roll the die when there are 3 sticks on the lid. Actually the optimum strategy is to stop rolling the die at the situations
(3/4/6), (3/5/5) and (3/6/4) but to continue to roll the die at all other situations with 3 sticks on the lid. The optimum strategy is denoted as\\ 

$\hat{S}_{13} = (1111111111/1111111111/111000111/00000000/0000000).$\\

Also for 14 and 15 sticks in the game with 3 sticks on the lid, the optimum strategy depends how the sticks are distributed amongst the two players. The optimum strategies are:\\

$\hat{S}_{14} = (11111111111/11111111111/1100000011/000000000/00000000)$ and 

$\hat{S}_{15} = (111111111111/111111111111/11000000011/0000000000/000000000).$\\

\begin{figure}[ht]
  \centering
  \includegraphics[width=0.75\textwidth]{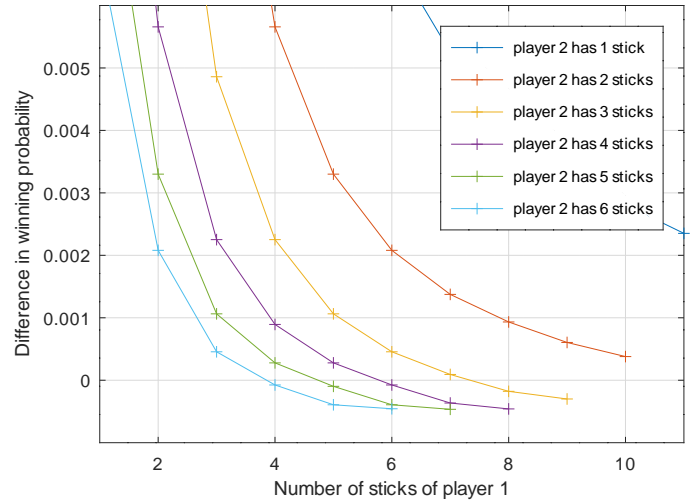}
  \caption{Difference in the probability to win when the player continues to roll the die or when they stop. There are 3 sticks on the lid, the number of sticks of the first player is given on the x-axis and the number of sticks of the opponent is a parameter indicated in different colors. If the difference is negative, it is better to stop.}
  \label{stop}
\end{figure}

The calculation of $\hat{S}_{15}$ took more than 24 hours on my personal computer and for 16 sticks it ran out of memory. Actually with so many sticks in the game the difference in the probability to win when you continue or when you stop and when there are 3 sticks on the lid is mostly below $5\cdot 10^{-4}$. Figure~\ref{stop} shows this difference as function of sticks of both players. The figure shows again that it is mostly better to continue when there are 3 sticks on the lid. It is only better to stop
\begin{itemize}
    \item when player 2 has 3 sticks (yellow or third curve from the right) and player~1 has 8 or more sticks, 
    \item when player 2 has 4 sticks (purple or fourth curve from the right) and player~1 has 6 or more sticks,
    \item when player 2 has 5 sticks (green or second curve from the left) and player~1 has 5 or more sticks,
    \item when player 2 has 6 sticks (blue or left curve) and player~1 has 4 or more sticks.
\end{itemize}

\section{Conclusions}

The strategies to maximize the winning chances in the game "Super Six" were determined for up to 15 sticks in the game. If there are 1 or 2 sticks on the lid, it is always best to continue to roll the die. If there are 5 sticks on the lid it is best to stop. If there are 4 sticks on the lid it is best to stop unless both players have one stick in their hands. If there are 3 sticks on the lid and not too many sticks in the game it is best to continue to roll the die, but with 13 or more sticks in the game it becomes advisable to stop depending on the situations listed in Chapter 8. 

\bibliography{references} 
\bibliographystyle{plain} 

\end{document}